\documentclass[12pt]{article}%
\usepackage{amsmath}
\usepackage{url}
\usepackage{amsfonts}
\usepackage{amssymb}
\usepackage{graphicx}
\usepackage{color}
\usepackage{url}
\usepackage{hyperref}%
\providecommand{\U}[1]{\protect\rule{.1in}{.1in}}
\providecommand{\U}[1]{\protect\rule{.1in}{.1in}}

\begin{document}

\title{From lecture to active learning: \\Rewards for all, and is it really so difficult?}
\author{David Pengelley}
\maketitle

In the centuries-honored \emph{I-You} paradigm, a mathematics instructor
provides first contact with new material via lecture, then expects students to
solve problems outside the classroom: \emph{I} lecture, then \emph{You} do
homework alone. Twenty-five years ago I began a personal journey away from
this paradigm, leading to a pedagogy that today fits within what is broadly
termed \emph{active learning}. I will provide an analysis of my  evolution to
a  particular non-lecture teaching philosophy. My experiences provide
encouragement for readers who may worry that alternatives to lecture are
complicated and time consuming. The main message is that it needn't be
difficult to create active learning for your students, and that there are
tremendous rewards for the instructor as well as for students. My goal is to
entice hesitant readers to take a teaching plunge.

Through the years my experiences convinced me that what can happen in an
active learning classroom can be greatly enhanced by good student preparation
before class. I\ expect students to prepare via reading, writing, and problem
work. Classroom activity can then build directly on their preparation, with
plenty of in-class feedback from fellow students and the instructor. Together
these components---student preparation and active classroom learning---enable
each student post-class to tackle higher-level homework. I call this
integrated approach my \textquotedblleft ABC method\textquotedblright, a
paradigm that I have refined in sixteen different courses at all undergraduate
and graduate levels. Lots of assignment examples, detailed guidelines for
students and for grading and daily logistics, and a holistic rubric are
available at \cite{pengelley}.

In a nutshell, to be described in detail later, my evolved paradigm for
student work is:

\begin{itemize}
\item A, \emph{due well before class: }Read, write questions, respond to my
questions. Graded for completion only.

\item B, \emph{bring to class}: Prepare `warm-up problems'. Graded for
completion in advance only.

\item \emph{During class}: Directly build active learning on both A and B.

\begin{itemize}
\item Lead brief discussion of the questions from A.

\item Compare and complete problems from B in groups: I facilitate, students present.
\end{itemize}

\item C, \emph{after class}: Complete a very few harder `final problems'
building on B from class. Marked carefully, these may be individually redone,
solely at my request, aiming for perfection. Holistic letter grade only. Final
level of achievement.

\item Parts ABC should comprise a very large part of the course grade; this
creates harmony between learning and evaluation, and reduces exams.
\end{itemize}

\noindent I will share my  thinking that led to this paradigm, its nuts and
bolts, what has and hasn't worked in which courses, how issues of time and
coverage work out, and student actions and reactions.

Equally important, since active learning likely seems shaky ground for those
primarily used to lecturing, I will provide some reassurance about how
demanding and time-consuming is the shift I\ made, and what the pitfalls are.
I will also address inertia, challenges, efficacy, teaching materials,
burnout, buy-in, and rewards, not only for students, but at least as
critically important and motivating, for instructors, since rewards for
instructors are perhaps crucial to overcome hesitancy.

In the intervening decades since my own evolution began, an enormous body of
research has developed that concludes  there are better alternatives than the
\emph{I-You} pedagogy for student success, collectively termed active
learning. Moreover, active learning confers disproportionate benefits for STEM
students from disadvantaged backgrounds and for female students in
male-dominated fields. And these benefits accrue while not disfavoring
high-achieving or more experienced students, or any demographic group
\cite{braun et al,CCI,laursen-2013,wieman-PNAS}. Recently the presidents of
fourteen professional mathematics societies joined to exhort us to shift from
\emph{I-You} toward active learning \cite{CBMS}. The question is how, and is
that hard to do?

All active learning paradigms\footnote{E.g., inquiry-based learning
\cite{laursen-2014,laursen-report}, interactive engagement \cite{CCI}, total
quality management \cite{williams}, just-in-time teaching \cite{jitt,novak},
peer instruction \cite{mazur}, flipped or inverted classroom \cite{flipped},
process oriented guided inquiry learning \cite{pogil}.} share two in-classroom
features. First, they reduce or eliminate lecture. Second, they devote
substantial classroom time to student involvement in mathematical work that
receives immediate feedback from other students and from the instructor. These
involve more of the ingredients \emph{You} and \emph{We}, and considerably
less of \emph{I}. Concomitantly, students will be more in charge of and
responsible for their own learning, while instructors will have increased
responsibility to guide student work. Within these broad parameters, many
variations are possible, so this is an exciting time of great experimentation
by many people who are seeking out and comparing a variety of good active
learning teaching techniques for mathematics. We already have a wonderful
resource in \cite{IPGuide}.

What I can contribute is one independently-developed example, the evolutionary
process of my own questioning, teaching experiences, reflection, and
adaptation leading to a philosophy and a longterm implementation within the
broad framework of active learning today. The distinguishing feature of my
evolved paradigm is its particular emphasis on tightly integrating pre-class
preparation via reading/writing and problem work with in-class active learning
and post-class follow-on homework.

\subsection*{Lecture}

From my early life as an \emph{I-You} student, I remember occasional
inspiration from lectures, but there was not much learning there that enabled
me to complete anything but rote homework. After all, lecture usually
primarily involves the instructor demonstrating that s/he can do the
mathematics. But my teaching experiences have led me to believe that this
rarely helps a student actually be able to do much mathematics, any more than
a swimming instructor demonstrating an hour of beautiful swimming techniques
successfully teaches a beginner how to swim various strokes. As a student, I
survived and prospered despite a lecture setting, but only by reading text
material repeatedly and integrating this with tackling homework challenges. I
now realize this was essentially autodidactical, my instructor's role chiefly
being to provide a schedule, expectations, homework feedback, and evaluation
via exams.

My subsequent decades teaching thousands of students suggests that few
students will very successfully self-teach in this way. The paradox for
readers of this article is that we are probably the most notable group of
exceptions; we are among the rare survivors or \textquotedblleft
thrivers\textquotedblright\ of the \emph{I-You} approach. But I expect we all
have frequent conversations with random adults, and with colleagues from other
disciplines, all former \emph{I-You} students of mathematics, in which we
receive strong unsolicited confirmation that the average \emph{I-You
}experience was a dramatic failure that left many scars.

During the years I lectured, many students told me \textquotedblleft I know
the math. I understand perfectly when you lecture, but then I can't solve
problems at home.\textquotedblright\ Of course in actuality this meant they
didn't really \textquotedblleft know the math\textquotedblright, but I didn't
know what I could do to help, other than to lead them through homework
problems. In retrospect, for all but possibly inspiration or rote learning, my
lecturing was ineffective, despite all my best efforts, and notwithstanding my
students' encouraging lauding of my lectures, their desire for them, and
belief in them. And since it wasted precious classroom time, it was
inefficient as well. In fact, classroom lecture may well become largely
obsolete, since with modern technology any recorded lecture can be viewed by
anyone, anytime, anywhere. How long will it take university administrators to
conclude that they need not employ professors to add more lectures to the
increasing number already archived? In short, professors had better have
something more to offer students than yet more lectures on settled subjects
\cite{bressoud-worst,lambert-mazur}. Of course we will all claim that our
students really do need much more than a lecture to succeed, and that we can
provide that. So isn't that what we should home in on? How then do we both
challenge students and guide and support their work as learners in truly
productive ways?

\subsection*{First contact with new ideas}

In rethinking the \emph{I-You} paradigm, much revolves around the question of
\textquotedblleft first contact\textquotedblright\footnote{I was enormously
influenced by the ideas of Barbara Walvoord \cite[pp. 53--63]{walvoord} on
first contact with new material, as beautifully described to me by Virginia
(Ginger) Warfield.}: How and when should a student first be exposed to new
material? In mathematics especially, absorbing and making sense of  new ideas
with any depth is usually a slow, highly individualized, intellectually messy
business. Lecture is by nature time-limited, one-size-fits-all, and totally
incompatible with the need to \textquotedblleft Stop, wait a minute, let me
think that through and pose a question.\textquotedblright\ In short, lecture
is on its face a poor means for first contact with demanding new material,
despite our natural inclination to the contrary, that as instructor we can
help students get started digesting new ideas by offering them a lecture
first.\footnote{For more on the disconnect between the role of lecture for
instructors and students, see \cite{weber}.}

So if in-class lecture provides poor first contact, then perhaps first
contact, and maybe even first problem work, might better occur \emph{before
}class, and something entirely different can happen \emph{during} class. This
could lead to the recently-named \textquotedblleft flipped\textquotedblright%
\ or \textquotedblleft inverted\textquotedblright\ classroom. In its original
conception, lecture and homework switch venues, with students watching
recorded lectures before class and working together on mathematics during
class. But watching recorded lectures has several drawbacks too. It is in many
ways at least as passive as watching a lecture live. And it suffers the same
cognitive drawbacks given above for in-class lectures, unless a student were
to frequently hit pause and replay, all the while thinking things through
deeply and asking questions, which is unlikely to happen without additional
pedagogical structure. 

My own conclusion, arrived at twenty-five years ago before recorded lectures
were even practical, was to evict lectures entirely and evolve new paradigms
instead. This shift was  catalyzed by my NSF-funded team's development of
student calculus projects \cite{srpc}. When we decided to allocate some time
for student groups to work in the classroom on their projects, with instructor
help, lecturing had to be reduced, and I began to expect students to read new
material before class instead. It is fascinating to realize that my shift
might never have occurred absent this external force.

\subsection*{Reading and the lecture-textbook trap}

How then do I want students to obtain meaningful first contact with new
material before class? My simplest answer in most courses is for students to
thoughtfully engage high-quality reading. And yet, while reflection and
thinking stimulated by reading can be extremely powerful, simply exhorting
students to read the book\ before class rarely works, since they seldom read
as suggested.

There is a gaping trap here, a truly vicious cycle in which students don't
read beforehand when they know the instructor will lecture, and instructors
lecture in large part because they know students haven't read. Breaking out of
this lecture-textbook trap was the most difficult teaching problem I\ ever had
to solve, but all else flowed from it. I felt it was my responsibility to
break this cycle by insisting (to self and students) that I will not lecture,
and instead arranging for in-class activity to be built on a foundation of
high-quality student preparation. A guiding motto was born: \textquotedblleft
Never lecture on something students can read instead.\textquotedblright

\subsection*{Written response to reading:\ Part A}

Resolving the lecture-textbook trap was the  catalyst for my entire journey,
enabling me over a number of years to evolve completely away from lecture.

An obvious resolution to the trap was to somehow convince students to read in
advance. However, seemingly making the trap even worse, reading alone is
insufficient, since they won't get much out of merely reading. I realized I
wanted students to reflect and think critically about what they read, to make
connections, and to respond in writing well before class. First, the nature of
an intellectually challenging writing task should hone their thoughtful
engagement and critical thinking and analysis, thus making the reading
worthwhile for learning. Second, I really need their written responses to
reading in order to prepare myself for the next productive non-lecture class session.

I have never asked students to summarize written material. Rather, I have
always challenged them to engage deeply by coming up with thoughts of their
own. For lower-division students, I find that written response to a couple of
well-crafted reading questions from me is essential. The questions I pose are
designed to stimulate students to read and think carefully, and to catalyze
and help guide class discussion. And sometimes my reading questions aren't
questions at all, but brief tasks based directly on the reading. I also expect
students at every level always to write their own good mathematical questions
about their reading, and additionally to write which new concepts are
confusing, what was well explained and interesting, what they had to reread
but eventually understood, and what connections they see to other ideas.
Textbook reading tends to provide polished answers to questions not even
meaningfully asked, and I attempt to get beyond this by expecting regular good
mathematical questions from students.

As examples of reading questions I might pose, in a first calculus
course,\ after reading an introduction to the derivative: \textquotedblleft
What are the different mathematical and physical interpretations we know of
for the derivative of a function?\textquotedblright\ and \textquotedblleft
Explain in your own words what your understanding is of the idea of the
derivative of a function.\textquotedblright\ Or in a discrete mathematics and
introduction to proofs course: \textquotedblleft Make up two great examples of
your own of multiply quantified statements, in which the meaning changes
dramatically when the order of the quantifiers is changed as in Examples 2.2.1
and 2.2.2. Explain why this is the case for each.\textquotedblright\ and
\textquotedblleft Make a good example of your own of each of the two types
(existential and universal) of multiply quantified statements discussed, and
then write and explain their negations.\textquotedblright\ Or for the
beginnings of a logic project based on primary historical sources
\cite{teaching-with-original-sources}: \textquotedblleft What is
logic?\textquotedblright, \textquotedblleft What did Boole attempt to
create?\textquotedblright, \textquotedblleft What is an
`implication'?\textquotedblright, \textquotedblleft How are implications
related to modern computers?\textquotedblright, \textquotedblleft According to
Aristotle, what is the difference between a sentence and a
proposition?\textquotedblright, \textquotedblleft What is a
syllogism?\textquotedblright\ More examples of reading questions are available
at \cite{pengelley}.

I mark each reading/writing assignment very quickly, holistically, with a
single +,$\checkmark$,- grade, only for seriousness of effort. I\ make as many
or as few comments as I wish or have time for, requiring only about five to
fifteen seconds per paper, since I am never reading detailed mathematics. My
greatest marking intent is to make sure that each student sees that I have
read and thought about what they wrote. Students become very faithful to this
reading and writing, and although I expect less than half a page of response,
some students become so emphatic about its benefits that they insist on
writing more, whether I want it or not! Some even explicitly credit their
success in the course to this activity.

If I receive student written reading responses up to one class period
beforehand, on paper or electronically, I can read them and determine how my
students are reacting to the new material. This best prepares me to guide
class without even a nagging impulse to lecture. I spend no time preparing a
lecture, rather I prepare notes on their writing so that I can best guide
their learning in the classroom.

Does this require reading material different from a textbook? Not necessarily,
provided the reading is genuinely accessible, interesting, stimulates
provocative thinking and questions about new ideas, and provides good grist
for class discussion. So I choose reading materials for these goals, possibly
utilizing multiple materials with different points of view to compare. This
does not mean that I choose material that promises to make the subject
\textquotedblleft easy\textquotedblright\ or a \textquotedblleft straight
path\textquotedblright, since such features often mean that the challenges,
questions, interest, and depth are missing, which does not serve learning in
the longterm.

\subsection*{In-class discussion of reading/writing: Building on Part A}

Class can now begin with a discussion directed by me, based on the few notes I
made while reading students' responses, which I\ first return with any
comments. It is always focused just on their writing, instead of a
shoot-in-the-dark lecture trying to address everything, without knowing what
students are struggling to understand. It can be geared specifically to meet
their needs for understanding the content I asked them to read, and the
discussion leaves no need for me to even consider lecturing. I sometimes have
students read out loud selected questions of theirs, and ask for class
reaction. This keeps the discussion focus on them and their thoughts, not
largely on mine. Since students have thoughtfully engaged the reading, this
second-contact in-class discussion never needs to be lengthy, usually 5-15\%
of class time, and the vast majority of time is available for something else.
What now could most usefully happen in class?

\subsection*{Warm-up problems beforehand:\ Part B}

Auspiciously, student written response to reading has prepared them for
productive initial mathematical problem work. So why not assign
easy-to-medium-difficulty \textquotedblleft warm-up\textquotedblright%
\ problems, also to be prepared in advance, and brought to class? In class
these problems can be compared, discussed, presented, and completed, using the
vast majority of classroom time, so that by the end of class the level of
student mathematical accomplishment, and their confidence, is high. I imagine
a traditional homework assignment, and put all but a few hardest problems into
the warm-up problems to be prepared for class. In terms of amount, it should
be just enough to keep everyone at work for the full class period, but not too
much more, since I want everyone to feel that this assignment is
satisfactorily completed by the end of class, and they are ready for a final
homework assignment of the harder problems.

This has been so successful that I have never had a single student express
reluctance about doing this homework as preparation before, rather than after,
class. The learning benefits quickly become obvious to students, since by the
end of class time, they are confident they have solved the easy-to-medium
warm-up homework, and feel ready to tackle harder problems at home.

\subsection*{In-class active learning: Integrating and building on Parts A and
B}

The stage has been set for classroom active learning by student preparation,
consisting of both reading/writing and warm-up problems. When I arrive at the
appointed hour, I find most students already comparing their prepared work in
groups of two to four. We begin together with a brief discussion based on
their reading/writing, as described earlier. Then much of class time is spent
in small group work refining their warm-up solutions, as I continually circulate.

My aim is to interact personally with every student or group multiple times.
My own classes have ranged from 10 to 50 students, with no lecture/recitation
format. Even in a class of 50, I am usually able to interact personally and
substantively with every student at least once during a 75-minute class period
that meets twice per week. I keep on the move, staying with each student or
group just long enough to provide encouragement, a little advice, and to learn
what they are struggling with. If a student or group is stuck, I help on that
point, then let them continue by themselves rather than rely on me for further
straightforward progress. If not stuck, then I may offer encouragement but not
stay long, unless a significant opportunity for enrichment presents itself to
be seized. I am always thinking about what I should do next. Should I select a
student to put a certain problem on the blackboard, or should I initiate a
whole-class discussion on a particular problem, or go on to another student or group?

I spontaneously initiate either whole-class discussions on particular
problems, or individual student board presentations and discussion. Often I
will ask several students to write solutions to various problems on the board
simultaneously, and then we discuss them all at once as a class. Sometimes a
writer is asked to verbally explain what has been written, sometimes not. Not
every problem gets presented or discussed.

I discovered painfully in one calculus course what happens if I steal
students' in-class work time by lapsing into lecture. I thought that the
material for the day was particularly tough, and that if I began with a bit of
lecture, it would help. After a while\ I saw frustration on my students'
faces, and I realized my mistake: They wanted to get to work on what they had
prepared for their valuable in-class time together, not listen to me. I now
realize I should have been very happy; they were in charge of their learning,
and they knew it.

I am generally laissez-faire about how groups arrange themselves. Three issues
with student group work are of concern. If the rare individual prefers to work
alone, I repeatedly encourage them to work with others, but I will ultimately
not force it. And if a group is not functioning optimally, I will sometimes
ask a member to switch to another group, to see if the dynamic is better.
Finally, I\ may occasionally reshuffle all the groups.

What are the forces ensuring that students really prepare problems before
class? Partly it is group peer pressure, and subtle pressure by me circulating
to observe student prepared work, and also the certain knowledge that I may
ask any student to present a prepared problem on the board at any time. I make
my presentation choices spontaneously, but very consciously, including to
press for better preparation beforehand by individuals if necessary by putting
them on the spot to present. But students' main motivation to prepare is their
experience that it creates a very effective learning environment, one in which
they will end class well equipped for the final, harder, after-class homework.

The warm-up problems are collected at the end of class, and marked
holistically +,$\checkmark$,-, again strictly for seriousness of effort at
preparation in advance, and they are important in the grade. (I could
alternatively have an extra copy due at the beginning or before class, e.g.,
photograph/scan and submit online.) Since the warm-up problems are dissected
in class, I never read them individually. I am interested solely in whether
the student prepared them in good faith before class. This literally takes
only five seconds per paper. Even though I am collecting them only at the end
of class, it is not hard to train myself to instantly assess preparation in
advance. This is particularly easy if, as often happens, there is a warm-up
problem that we didn't get to in class discussion; then I can easily see on
each student's paper whether the problem was prepared beforehand or not.

\subsection*{Final problems after class:\ Part C, exams, and course grade}

With student preparation before class of reading/writing and warm-up problems,
all in support of in-class discussion, group work, and presentation, it
remains only to assign a very few (two or three) harder \textquotedblleft
final\textquotedblright\ homework problems for students to complete after
class. These are like the hardest few of traditional homework, but now build
on what students have already achieved before and during class. The final
problems\ are the only daily work needing detailed marking, representing each
student's highest level of achievement on the day's material.

Their papers receive prompt and very careful feedback, a single holistic
letter grade\footnote{I have long found using points in marking to be a
time-wasting, exhausting, distracting, and deceptive morass that sends the
wrong message to students and invites trouble. While I may write a lot on
student papers, I always assign only a single holistic qualitative evaluation
to a paper, be it homework or an exam.}, and possible prompt redoing of
individual problems at my initiative to bring to perfection. They are normally
never discussed in class. These higher-level problems are at the core of a
student's course grade; I consider them the best measure of what each student
has learned and accomplished. The message to students is that their three
daily written components are the fundament of both learning and evaluation, so
I find it critical that they form the vast majority of the course grade. This
leads to a reduction in exams, which I believe is good, because my experience
is that timed in-class exams are a poor way for most students to demonstrate
what they can actually do.

I have always made the three ABC components together of student daily work at
least 60\% of the holistically-assigned course grade, with the carefully
graded final problems\ (Part C) dominant in my mind. Since almost all of the
reading/writing (A) and problem preparation (B) assignments earn a +, the
harder, after-class, final problems (C) with letter grades tend to become
paramount in each student's course grade. In some courses I still use some
form of midterm and/or final exam, while in others, especially upper division
courses, I sometimes use no exams, since the daily Part C graded work with
high expectations is more than sufficient for evaluation. Usually I leave open
the option for exams, and often decide against them as the term progresses.

A word of warning based on experience: Once upon a time I didn't clearly
separate the warm-up from the final problems, but this led to complications,
including lesser student effort on the warm-up problems before class; I find
it works far better with the two sets well separated.

How can one be sure that these final problems completed outside class
represent the work of each individual, especially since I encourage students
to work together in class so much? Truth be told, even on the final problems I
make it clear that students may discuss them together. My ironclad rule,
though, is that when they go to write them up, this must be done alone, so
that no two papers should look alike. Since these harder problems are never
just a calculation or formulaic, but always involve explanation of ideas, I
can detect not only the level of understanding, but also easily observe if two
papers are too similar. If so, which occasionally happens at the beginning of
the course, I speak with the students involved to reiterate my expectations as
emphatically as necessary. Another issue could be students finding and copying
solutions from elsewhere, e.g., online, to problems for Parts B or C. Since
this completely defeats their purpose, if this is a danger one must create or
tweak problems to avoid it.

\subsection*{Pr\'{e}cis of student assignments}

To sum up the evolved paradigm, which tightly integrates before-, during-, and
after-class work, students write three homework papers for each daily unit of
content (in my case twice weekly), which I call parts A,B,C, and which replace
the \emph{I-You} paradigm:

\begin{itemize}
\item \emph{You}, Part A: Read/write, received up to one class period early
for me to prepare for leading class discussion. Marked quickly +,$\checkmark
$,- for effort only.

\item \emph{You}, Part B: Warm-up problems, prepare and bring to class. Marked
quickly and holistically +,$\checkmark$,- for preparation only. Submit in
person, and/or online photo/scan before class.

\item \emph{We}: In-class discussion, group work, presentations, all built on
Parts A,B for the given unit.

\item \emph{You}, Part C: A very few harder final problems, completed after
class and written up alone. Marked carefully with feedback and holistic letter
grade, sometimes specific problems redone at my request.

\item Together Parts A,B,C constitute the majority of the course grade,
reducing or eliminating exams.
\end{itemize}

Note that for a given unit, Parts A,B,C are due at different times, so there
is a rolling nature to the coverage of multiple units. Students easily adapt
to this provided I give clear guidance, and it has integrative benefit.

\subsection*{Inertia}

There are many forces binding an instructor to the \emph{I-You} paradigm, even
if she is open to change. First, we naturally tend to teach as we were taught.
Second, a lecture \textquotedblleft covering material\textquotedblright\ is in
the tradition of fulfilling professorial duty. It is hard to let these go, and
to realize that instructor coverage\ does not necessarily help students do
mathematics. Students, too, are mostly happily complicit, generally very
comfortable with passive receipt of a lecture. Certainly it is much easier
than having to do any actual work in the classroom, and they can believe they
must have learned something from lecture.

Third, it takes real effort to change pedagogy, and the change will likely
catalyze a process of further evolution, so it is a major commitment. Fourth,
there is an element of uncertainty, worry, and perhaps fear of classroom
disaster. Lecture is well-known, often easy, predictable, and contains
essentially no element of risk, since it is totally controlled by the
instructor based on preparation in advance, with little chance of something
unexpected or surprising from students. On the other hand, in a classroom of
continual interaction with students, reacting to, adjusting, and guiding what
students initiate may seem scary. Moving away from lecturing amounts to
relinquishing total control, but hopefully without totally losing control,
since one still has overall guiding responsibility. Creating the right balance
is a challenge.

Therefore, since shifting from \emph{I-You} requires overcoming much inertia,
it will most likely occur only if one sees large benefits and rewards (and not
too many scary challenges) for both student learning and for instructors
themselves; here I would like to offer much encouragement from experience.
Let's begin with the students.

\subsection*{Benefits for students and learning}

My experience is that students respond well to preparatory work provided they
quickly and consistently experience the advantages, know it is highly valued,
i.e., in class and in their grade, and are fully expected to contribute in
class based on their preparation in advance. They then find their in-class
work time valuable, engaging, rewarding, often exciting, and
confidence-building. Completing the warm-up problems with feedback from me and
fellow students in the classroom prepares them well for success with the few
final, harder problems to be completed after class for careful grading, and
they know and greatly value that. Many times my students are so absorbed in
group work completing warm-up problems that they don't realize when class time
has ended. I apologetically interrupt the whole class to tell them that class
ended five minutes ago! When does that ever happen with a lecture?

The reduction of exams along with the predominant emphasis on daily work for
both learning and course grade creates a much steadier workload for students,
yielding the cognitive advantages of spaced learning, and relief from the
typical cram/exam/forget phenomenon that doesn't foster longterm learning.
This also places learning and evaluation in harmony, reducing stress and
producing more consistent quality of work. My impression from a many-years
evolution is that, with these approaches, my students work more, and more
successfully learn course material.

Student course evaluation comments are almost uniformly positive about the
pedagogy, and indicate a high level of buy-in. They typically credit
preparation in advance for in-class collaborative work as extremely effective
for their learning, and for keeping them on top of the course with less
stress. Students also often remark that the emphasis on student participation
makes the subject come alive. Quite frequently they ask why other mathematics
courses are not taught this way.

I mention here one anecdote that still astonishes me, from an abstract algebra
course intended both for mathematics majors and future secondary mathematics
teachers. Although the entire course was focused on mathematics, at the end of
the semester one student came to my office to tell me that for her, more
important than the mathematics had been the teaching style, and that she had
consciously spent the entire term studying the pedagogy, with the aim of
adapting it in her own teaching. Never had I dreamt that while thinking I was
teaching abstract algebra I was actually also teaching pedagogy.

\subsection*{Challenges for instructors}

Shifting from \emph{I-You} to something like \emph{You-You-We-You} has initial
challenges for an instructor. As with anything new, more effort will be needed
the first time. Experience pays off handsomely, though, and after once or
twice through, my experience is that the overall workload should be no greater
than for \emph{I-You}.

It is critical that students have confidence from the start. I build this by
briefly explaining to the class the evidence for how an active learning
paradigm will enable them to be successful in the course, that class time will
be interesting, productive, and satisfying, that it will prepare them well for
the harder homework, and that this daily work is the great majority of their
course grade. And I assure students that I will be there to give personal help
in class every day. Then I watch and listen to how things are going,
especially in the first weeks, repeat my explanation for active learning as
necessary, and take prompt steps to resolve any confusion and alleviate discomfort.

I have benefited from the fact that I match students' stereotypes of who is a
professor and an authority in ways that many other instructors will not.
Implicit bias affects students' perceptions of their teachers' expertise, and
thus their willingness to extend the benefit of the doubt to a class structure
that does not conform to their prior notions. So some instructors will need to
be even smarter about strategies for obtaining student confidence and buy-in
and managing resistance. My principal advice is to keep student confidence
always in the forefront of one's mind.

One must learn, as addressed above under In-class Active Learning, how
mindfully to make decisions that support student learning in a less
predictable classroom environment where control and responsibility is being
loosened and partially handed to students. I also need to keep reminding
myself that in a nonlecture classroom, it is students who should be doing the
mathematics, not the instructor, since I already know the mathematics, and
they are the workers and learners. My job really should be that of effective,
efficient, encouraging, and hopefully inspiring, guide and
manager.\footnote{Frank Williams was a seminal influence for me, for which I
am forever grateful. It was he who helped me understand the proper role of the
instructor to complement the student's role of worker and learner
\cite{williams}.} Neither should anyone expect their learning to be easy: I
can be helpful in many ways, but the learning is their work, just as when
Euclid is said to have replied to King Ptolemy's request for an easier way of
learning mathematics that \textquotedblleft there is no royal road to
geometry\textquotedblright.

Perhaps the greatest danger for an instructor is that with students handing in
homework Parts A,B,C for each class day, it would be all too easy for me to do
way more homework marking than I should, and therefore spend more time
teaching this way; I\ have witnessed colleagues insistently fall into such a
hole when trying this approach. Each of the three parts is crucial for student
learning, but Parts A and B do not need grading or instructor feedback on the
mathematics, since this all happens in class. While Part C is carefully
marked, and perhaps pieces redone, it consists only of two or three harder
problems, making grading manageable.

I am often asked how to start on the first day of a term, since Parts A and B
are to be completed before class. First, I never lecture; instead I model the
pedagogy on that first day by having students work together on meaningful
mathematical activity. Then between the first and second class days, I have
students submit their first response to reading, the only time anything
happens off schedule; thus we are ready for a normal routine on the second day
of class, when the first warm-up problems and the second response to reading
are due.

Is my paradigm a one-size-fits-all approach? While I have found the basic
components to be universally successful, the details may best differ between
courses at different levels, or with different meeting schedules, or with
large-size classes.

For instance, in a mathematics appreciation general education course at the
lowest college level, I emphasize hands-on activity more than reading, for
both work at home and in class; and after-class work often entails students
just writing up what they discovered in class. At the other extreme, in a
Ph.D. level graduate course, I often ask students to contrast multiple
different written approaches, and in class I will ask students to present
their own versions of proofs at the board and lead discussion thereof.

Regarding class meeting schedules, in the past several years I\ have worked
individually with numerous faculty and graduate teaching assistants who are
adapting this pedagogy to different schedules than mine, which has always
consisted of two 75-minute class sessions weekly. I see that the scheduling of
Parts A,B,C may best be adjusted, with some consolidation in a course with
three or more meetings per week, e.g., only two Parts A and B weekly, and/or
one Part C. Some adaptations are very innovative and substantial, such as in
\cite{dunmyre}.

For large class sizes, or courses traditionally scheduled in
lecture/recitation format intended specifically for lecturing, I see other
successful adaptations being made, e.g., with one or multiple graduate
teaching assistants and/or undergraduate learning assistants
\cite{Learning-Assistant-Alliance,Learning-Assistant-Program} present in the
classroom with the instructor to work with students.

Amongst instructors adapting this pedagogy, the vast majority have had great
success, while a few have struggled. My perception is that those who struggled
had unintentionally combined old and new pedagogies in incompatible ways. For
instance, if daily in-class active group work is graded for correctness rather
than for participation, completion, and effort, this can undermine the
learning process by shifting student attention entirely towards a correct answer.

Finally, what must I at minimum have within my control in order to teach this
way? I need my students to have access to good reading and problem material
that I can assign as needed, including reading/writing and problem preparation
in advance of class. I need the daily pre- and post-class assignments to be
the core of students' work and grade. And I need to be able to mold the
classroom environment into an active one and gain the confidence of my
students. All else is flexible.

\subsection*{Benefits and rewards for instructors}

Perhaps for many instructors, at the end of the day it will also be the
personal rewards, not just those for students, that will seem attractive about
alternatives to \emph{I-You}. I admit I have reaped tremendous personal rewards.

Class time has higher-quality interactions and is more exciting when one is
frequently discussing interesting mathematical ideas with individuals and
groups, and they are coming up with questions and ideas and points of view
that one hadn't anticipated. The enthusiastic response from students is
extremely gratifying, as is the learning success one sees. In short, I\ enjoy
interacting with my students much more, a huge benefit!

Marking student work is more rewarding in two senses. First, exams are fewer.
Second, marking time is spent primarily on the few harder homework problems,
which are more interesting to mark, not on the easier material that has been
dealt with in class. And the remaining time is spent mostly reading student
responses to reading, which stimulates and prepares one with confidence to
lead a class discussion most useful for student learning.

Time, ah time: My experience over many courses is that an alternative to
\emph{I-You} need not take more instructor time overall, provided one does not
fall into the trap of unnecessary over-marking of student homework prepared
for class. A perhaps surprising timesaver is that students often need less of
my time in office hours: By replacing lecture with student interaction with
each other and with me on active work inside the classroom, students get most
of the help they need, and their questions answered, in class. Moreover, the
steadier workload mentioned above applies to instructors as well, so there is
very little end-of-term stress, and no longterm burnout.

With rewards as strong as these, I could never return to \emph{I-You}. Carpe diem!

\subsection*{Is there really an elephant in the classroom?}

Finally, consider the question of coverage, an intimidating and much-feared
elephant. When I talk with \emph{I-You} instructors about replacing lecture
with student work in class, they almost invariably reply \textquotedblleft But
then I couldn't cover all the material in the syllabus\textquotedblright. My
primary answer, of course, is that it is not the instructor who needs to cover
the material, but rather the student.

I have found, in teaching many types and levels of courses, that if
high-quality first contact and initial mathematical work happens before class,
thus making lecture irrelevant and redundant, and if class time is instead
used for student work with others and with the instructor to build on the work
prepared in advance, then coverage is always more efficient, not less so. To
me this simply makes logical sense: If lecture is a largely ineffective use of
precious classroom time for student learning, then offering students a guided
active-learning classroom environment, working with each other and with me,
seems likely to proceed more efficiently, especially when first-contact
reading and preparatory work happens before class. Specifically, I have taught
this way in first-year calculus courses with multiple sections all following
the same lockstep routine with common exams, where students in my section had
to progress at the same rate that other instructors were lecturing, and this
was no problem at all. In fact it was in exactly that setting, with a class of
45 students and no grader or teaching assistant, where I first developed and
refined the approach described here.

My consistent experiences after transforming \emph{I-You} into
\emph{You-You-We-You}, in many courses at all levels and for all college
audiences, is that the content is actually less rushed. I found no fearsome
coverage elephant in the classroom as I redesigned it, even with the same
syllabus as other instructors.

\bigskip

\bigskip

\noindent\textbf{Acknowledgements}: I am indebted to Sandra Laursen, Barbara
Walvoord, Virginia Warfield, and Frank Williams for critical insights along my
journey, and to Pat Penfield and many others for encouragement and
constructive criticism.

\bigskip

\noindent\textbf{David Pengelley} is professor emeritus at New Mexico State
University, and courtesy professor at Oregon State University. His research is
in algebraic topology and history of mathematics. He develops the pedagogies
of teaching with student projects and with primary historical sources, and
created a graduate course on the role of history in teaching mathematics. He
relies on student reading, writing, and mathematical preparation before class
to enable active student work to replace lecture. He has received the MAA's
Deborah and Franklin Tepper Haimo teaching award, loves backpacking and
wilderness, is active on environmental issues, and has become a fanatical
badminton player.

\bigskip

\noindent\textit{Department of Mathematical Sciences, New Mexico State
University, Las Cruces, NM 88003\newline Department of Mathematics, Oregon
State University, Corvallis, OR 97331\newline davidp@nmsu.edu}
\end{document}